

\def\proof{\noindent{\bf Proof.  }}

\def\O{{\cal O}}
\def\L{{\cal L}}
\def\F{{\cal F}}
\def\E{{\cal E}}
\def\oqy{\Omega^q_Y}

\def\oqqqy{\Omega^{q-1}_Y}
\def\oqx{\Omega^q_X}
\def\oqqx{\Omega^{q+1}_X}
\def\oqqqx{\Omega^{q-1}_X}
\def\oqyx{\Omega^q_{Y\vert X}}
\def\oqqyx{\Omega^{q+1}_{Y\vert X}}

\def\oqqqyx{\Omega^{q-1}_{Y\vert X}}
\def\Diag{\def\normalbaselines{\baselineskip25pt
\lineskip3pt\lineskiplimit3pt}
        \matrix}
\def\lto{\rightarrow}
\magnification=1200
\vglue.5in

\vskip .5 cm

{\centerline{{\bf On stability of tangent bundles of Fano manifolds
with $b_2=1$.}}}
\vskip .5 cm
{\centerline{Thomas Peternell and Jaros\l aw A.~Wi\'sniewski}}
\vskip .5 cm

\beginsection Introduction

An important outstanding problem in differential geometry
asks which Fano manifolds
$X$ with $b_2(X) = 1$ admit a K\"ahler---Einstein metric.
The only known general result (besides
some necessary conditions as the reductivity of
the Lie algebra of holomorphic vector fields and
the Futaki invariant) is
the positive answer for rational-homogeneous manifolds,
which holds also if $b_2(X) > 1.$ A weaker
and more algebraic question would ask
whether the tangent bundle $T_X$ is stable (with
respect to the anticanonical bundle $-K_X$).
For the connection between both notions see e.g.~[Ti].
 If $b_2(X) > 1$, $T_X$ is not necessarily
stable due to the geometry of contractions of extremal rays;
for 3-folds this has been studied
completely by A.~Steffens in his thesis [St].
There are several papers considering the 2-dimensional case,
i.e.~del Pezzo surfaces, see [Ti] for references.
This paper contributes to the stability problem
in case $b_2(X) = 1.$
\medskip
Of course the stability of $T_X$ is very much related to
cohomology vanishings of type
$$H^p(X,\Omega^q(t)) = 0.$$
These vanishings are studied in Section 1. We say that a pair
$(Y,{\cal O}(1))$ consisting of a projective manifold
$Y$ and an ample line bundle
${\cal O}(1)$ has special cohomology if

{\narrower\medskip
\item{(a)} $H^p(Y,\Omega_Y^q(t)) = 0 $ \
{\rm for} \  $0 < p < dim Y$, $p+q \ne dim Y$ and
$ t \in {\bf Z}\setminus \{0\}$
\item{(b)}$ H^p(Y,\Omega^q_Y ) = 0$ \  {\rm for} \
$ 0 < p < dim Y , p+q \ne dim Y, p \ne q$
\item{(c)} $H^p(Y,\Omega^p_Y) = {\bf C}$  \ {\rm for  } \  $ 0\leq p \leq dim
Y$,
\medskip}

In other words, conditions (b) and (c) say that the odd Betti numbers except
the
are 0 while the even Betti numbers are 1.

Of course $({\bf P}^n,{\cal O}(1))$
has special cohomology. But due to the following
``functorial'' theorem we have much more examples.

\proclaim Theorem 1. Let $(Y,{\cal O}(1))$ have special cohomology.
\item{(1)} Let $X \in \vert {\cal O}_Y(d)\vert $
be a smooth divisor for some $d > 0$.
Then $(X,{\cal O}(1)_{\vert X})$ has special cohomology.
\item{(2)} Let $\pi : X \longrightarrow Y $
be the $k$-cyclic cover branched
along some $D \in \vert {\cal O}(kd)\vert $.
 Then $(X,\pi^*({\cal O}_Y(1))$ has
special cohomology.
\item{}

 This theorem concerns vanishing of higher cohomology groups.
The same methods yield a
 $H^0$-vanishing theorem which can be stated as

\proclaim Theorem 2.  Let $(Y,{\cal O}(1))$ have special cohomology,
moreover assume the
tangent (or cotangent) bundle of $Y$ to be semi-stable.
\item{(1)} Every smooth $X \in \vert {\cal O}_Y(d))\vert $
has stable (co-)tangent bundle, if dim$X \geq 3$
\item{(2)} Let $\pi : X \longrightarrow Y $ be the
$k$-cyclic cover branched along some
$D \in \vert {\cal O}(kd) \vert $, $k \geq 2.$
Then $X$ has stable (co-)tangent bundle,  if dim$Y \geq 3$.
\item{}

Of course Theorem 2 does not solve the stability
problem for Fano manifolds completely.
However we can prove in Section 2:

\proclaim Theorem 3. Let $X$ be a Fano n-fold of index
$r$ with $b_2(X) =1$.
\item{(1)} If $r \geq n-1$, then $T_X$ is stable.
\item{(2)} If $r =  n-2$ and if there are ``enough'' smooth
divisors in
$\vert H \vert $, $H$ being the ample generator of
$Pic(X) = {\bf Z}$, then $T_X$ is stable .
\item{(3)} Every Fano 4-fold with $b_2(X) = 1$
has stable tangent bundle

The existence of enough smooth divisors
(for the precise statement see Section 2) in (2) is expected
to hold on every Fano $n$-fold of index $n-2$.

It should be mentioned that the methods of Section 2 are to
some extend ad hoc methods relying on classification results
of Fujita, Mukai and others, so
they do not apply immediately to manifolds of lower index.
\bigskip

We would like to thank Professors P.~Pragacz
and A.~Daszkiewicz for very
useful discussions and informations on the
cohomology of Grassmannians. We are also indebt to
Professor H.~Flenner for discussions concerning weighted
complete intersections as well as for informing us
about his results from [F]. Last but not least we would like to
thank Professor D. Snow for kindly informing us about the cohomology
of the $10-$dimensional spinor manifold and the referee for saving us from
an error in an earlier version.

Part of this research has been
granted by the Schwerpunkt ``Komplexe Mannigfaltigkeiten'' of
the Deutsche Forschungsgemeinschaft; the
second-named author would like to acknowledge a support
from a Polish grant KBN GR54.

\beginsection 0. Preliminaries

For a complex manifold $X$ we denote by $\Omega_X^q$
the sheaf of holomorphic q-forms; we will use often
 $\Omega_X$ instead of $\Omega_X^1.$
$T_X$ will always denote
the tangent bundle (sheaf) of $X$ and for a smooth subvariety
$Z\subset X$, $N_{Z/X}$ will denote the conormal bundle.
The restriction of a coherent sheaf ${\cal F}$ over $X$ to $Z$
will be denoted by ${\cal F}_{\vert Z}$.
We will frequently identify divisors on $X$ with line bundles
associated to them. For a line bundle ${\cal O}(1)$
and a coherent sheaf $\F$ on $X$ by $\F(t)$ we will denote
$\F\otimes{\cal O}(1)^{\otimes t}$.
\medskip
Stability or semi-stability is understood in the sense of
Mumford---Takemoto, so a vector bundle ($=$
locally free sheaf) $\E$ is (semi-)stable with respect
to an ample line bundle $L$ if for all
coherent sheaves $\F\subset \E$, $0 < rk\F < rk\E$, it holds
$$ \mu (\F) < \mu (\E) \ \ \
(\hbox{respectively, }\ \mu (F) \leq \mu (\E)),\
\hbox{ where } \mu (\F) = {{c_1 (\F).L^{n-1}}\over {rk\F}}.$$
If $b_2(X) = 1$ and the ample line bundle is not mentioned
explicitly we mean of course ``the'' ample
generator of $PicX$.
\medskip

A Fano manifold has by definition ample anticanonical
bundle $-K_X = det(T_X)$. Its index is the largest
positive integer $r$ such that ${1 \over r} K_X \in PicX$.
The coindex $c$ is given by
$c = dimX-r+1$.
By the Kobayashi-Ochiai theorem $c$ is always non-negative,
$c=0$ iff $X={\bf P}^n$,
and $c=1$ iff $X={\bf Q}^n$, the $n$-dimensional smooth quadric.
\bigskip

In Section 2 we will use informations
on various analytic cohomology groups of Grassmannians,
in particular for $G(1,5)$, the Grassmannian of lines in ${\bf P}^5$.

\proclaim Lemma 0.1. Let $G = G(1,n)$ be the Grassmann
variety of lines in projective $n$-space. For
a real number $s$ we denote by $\lceil s\rceil $ its round up,
by $\lbrack s\rbrack$ its round down. Let
${\cal O}(1)$  be the ample generator of $Pic(G)$.
Then $H^p(G,\Omega_G^q(t)) \ne 0$
if and only if one of the following conditions is satisfied
\item{(1)} $t=0$ and $p=q$
\item{(2)} $p=0$ and
$t\geq \hbox{min} \bigl( q+1,\lceil {q\over 2} \rceil +2 \bigr)$
\item{(2')} $p=dim G = 2(n-1) $ and $t\leq \hbox{max}
\bigl( -2n+q+1,\lbrack {q\over 2}\rbrack -n-1 \bigr)$
\item{(3)} $0 < p < 2(n-1)$, $p+1 \leq t \leq n-p$ and $q=2p+t-1$
\item{(3')} $0 < p < 2(n-1)$,
$-p-1\leq t \leq -n+p$ and $q=2p+t-2n+3$.

The proof of (0.1) can be deduced from [La, Cor.~4.1],
see also [We, 2.4] and [Sn].

\bigskip

In section 1 we will discuss some extensions of the following result
of Flenner.

\proclaim Theorem 0.2. {\rm [F, Satz 8.11]}
Let $X$ be a weighted complete intersection of dimension $n$
in a weighted projective space. Let ${\cal O}(1)$
denote the restriction to $X$ of the universal
${\cal O}(1)$-sheaf from the weighted projective space.
If $X$ is smooth then
\item{(A)} $H^q(X,\Omega^q_X)={\bf C}$ for $0\leq q\leq n$,
$2q\ne n$, and
\item{(B)} $H^p(X,\Omega^q_X(t))=0$ in the following situations:
\itemitem{(a)} $0<p<n$, $p+q\ne n$, $p\ne q$;
\itemitem{(b)} $0<p<n$, $p+q\ne n$, $t\ne 0$;
\itemitem{(c)} $p+q>n$, $t>q-p$;
\itemitem{(d)} $p+q<n$, $t<q-p$.
\medskip

The knowledge of the above cohomology allows to set the stability
of the tangent bundle, so that we have the following:

\proclaim Corollary 0.3. Let $X$ be as above. Assume moreover
that $n\geq 3$. Then $\Omega_X$ (as well as $T_X$) is stable.

The argument to prove the corollary is typical and will be
used several times throughout the paper:
For $n\geq 3$ we have $b_2(X)=1$ and thus we may assume
that $det\Omega_X\equiv {\cal O}(-r)$ and moreover
that $r\leq n-1$ because otherwise $X$ is either a projective
space or a quadric.
By $d$ det us denote the self-intersection $({\cal O}(1))^n$.
Let $\F\subset\Omega_X$ be a
subsheaf of $\Omega_X$ of rank $q$,
 $1\leq q<n=rk\Omega_X$
and $det\F\equiv{\cal O}(k)$.
Then
$$\mu(\Omega)=-{r\cdot d\over n}\ \ \hbox{ and }\ \
\mu(\F)={k\cdot d\over q}$$
hence we are to prove that
$$-k>q\cdot {r\over n}.$$
We may assume that $\F$ is reflexive so that
$det\F\subset \Omega^q_X$ is an invertible subsheaf.
Therefore we have a non-zero section
in $\Omega^q_X(-k)$ and from the previous theorem,
point {\sl (d)}, it follows that
$-k\geq q$. Since $r/n<1$ this proves the desired inequality.

\beginsection 1. Vanishing theorems on complete
intersections and cyclic coverings.

In the present section we discuss some results on vanishing
of cohomology of sheaves of twisted holomorphic differential forms
and their wedge products. To make citations
easier let us introduce the following

\proclaim Definition. A pair $(Y,\O(1))$ consisting of a projective
manifold $Y$ of dimension $\geq 3$
and an ample line bundle $\O(1)$ over $Y$
has special cohomology if the following conditions are satisfied
\item{(a)} for $0<p<dimY$,
$p+q\ne dimY$ and $t\in {\bf Z}\setminus \{ 0\}$
$$H^p(Y,\oqy(t))=0,$$
\item{(b)} for $0<p<dimY,  p+q \ne dimY ,  p\ne q$
$$ H^p(Y,\Omega^q_Y) =0$$
\item{(c)} for $0\leq p \leq dimY, 2p \ne dimY $
$$ H^p(Y,\Omega^p_Y ) = {\bf C}.$$
\item{}
\par

Note that, due to the well-known Bott formulae, the pair
$({\bf P}^n,\O(1))$ has special cohomology.
Flenner in [F, Satz 8.11] proved that smooth weighted
complete intersections in weighted projective
spaces have special cohomology, too.
Our point is to prove
that the class of manifolds having special cohomology contains
also complete intersections in them and cyclic coverings of them.

\proclaim Theorem 1.1. Let $Y$ be a projective manifold
of dimension $n+1$ and $\O(1)$ an ample line bundle
over $Y$. Assume that the pair $(Y,\O(1))$ has
special cohomology.
Let $X\subset Y$ be a smooth divisor from a linear system
$\vert {\cal O}(d)\vert$. Then the pair
$(X,\O_X(1))$ has also special cohomology.

\medskip

\noindent{\bf Remark.} If $p+q>n$ and $t>0$, or
$p+q<n$ and $t<0$, the vanishing in the condition (a) of the above definition
is
immediate from the Kodaira-Nakano vanishing theorem.
Also let us note that for $t=0$ and $p+q<n$ we have an isomorphism
$$H^p(Y,\oqy)\rightarrow H^p(X,\oqx )$$
(restricting of $(p,q)$-forms) which follows from
Lefschetz hyperplane section
theorem. Thus, the pair $(X,\O_X(1))$ as in the theorem satisfies
the conditions (b) and (c) from the
above definition and only (a) has to be proved outside the range which is
settle
Kodaira-Nakano.
\par
Also, let us note that the theorem is clearly true for
 $q=0$; this follows from the cohomology of the sequence
$$0\longrightarrow \O_Y(t-d)\longrightarrow
\O_Y(t)\longrightarrow \O_X(t)\longrightarrow 0.$$

\medskip

Before proving the theorem let us recall
two exact sequences of sheaves,
the first one on $Y$:
$$0\longrightarrow \oqy (t)\longrightarrow\oqy (t+d)\longrightarrow
\oqyx (t+d)\longrightarrow 0\eqno(1)$$
and the next one which is a twisted exterior power
of cotangent---conormal
sequence on $X$:
$$0\longrightarrow \oqx (t)\longrightarrow\oqqyx (t+d)\longrightarrow
\oqqx (t+d)\longrightarrow 0.\eqno(2)$$
We will use these sequences to express cupping with $c_1\O(d)$,
compare with [Mt, (2.2)] and also [So, (0.7)]:

\proclaim Lemma 1.2. The composition of
appropriate maps on cohomology
of the above sequences:
$$\eqalign{H^{p-1}(Y,\oqqqy)\rightarrow H^{p-1}(X,\oqqqyx )\rightarrow
H^{p-1}(X,\oqqqx )\cr \rightarrow H^{p-1}(X,\oqyx (d))\rightarrow
H^p(Y,\oqy)}\eqno(3)$$
is  cupping with $c_1({\cal O}(d))$ (and thus is an isomorphism
for $p+q<n+1$).
\par

\noindent{\bf Proof of lemma. }
The above sequences come by wedging those for $q=1$
so first we examine the situation for $q=1$.
Let us assume that the divisor $X$
is defined locally on $Y$ on a covering $(U_i)$
by functions $(f_i)$.
The morphism
$H^0(Y,\O)\cong H^0(X,\O)\lto H^0(X,\Omega_{Y\vert X}(d))$
in (2) is defined by differentials
$(df_i)$ which glue to a section over $X$.
Since the first map in (1) is locally defined by multiplying
by $f_i$, the boundary map in the cohomology of (1)
associates to the unit section the following 1-cocycle
$$({df_i\over f_i}-{df_j\over f_j})=(d(log({f_i\over f_j}))).$$
This, however, is the first Chern class
of the divisor $X$, see [H, exercise III.7.4]
and [Mt].

We check that this description extends for other $p$ and $q$:
for a \v Cech cocycle $\alpha\in\check Z^{p-1}((U_i),\oqqqy)$
the composition of first 3 arrows associates a cocycle
which on $U_{i_0}\cap\dots\cap U_{i_{p-1}}$ assumes value
$df_{i_0}\wedge\alpha_{i_0\dots i_{p-1}}$
(note that $df_{i_0}$ is restricted
to the intersection $U_{i_0}\cap\dots\cap U_{i_{p-1}}$).
The boundary map then gives the cocycle
$${df_{i_1}\over f_{i_1}}\wedge \alpha_{i_1\dots i_p} -
\sum_{k=1}^{k=p} (-1)^{k+1} {df_{i_0}\over f_{i_0}}
\wedge\alpha_{i_1\dots i_{k-1}i_{k+1}\dots i_p}$$
which, since $\alpha$ is a cocycle, is equal to
$$({df_{i_1}\over f_{i_1}}-{df_{i_0}\over f_{i_0}})
\wedge \alpha_{i_1\dots i_p}$$
and this is just the $\wedge$-product of the first Chern class of
$X$ with the cocycle $\alpha$.
\medskip

\noindent{\bf Remark.}
Alternatively, one may use the Dolbeault complex to prove the above
lemma.

\medskip
\noindent{\bf Proof of Theorem 1.1.}
First let us note that we may assume
$p+q<n$ because then the other case is set by Serre duality.

Since we know the vanishing of $H^p(X,\oqx (t))$ for $t<0$ and
$p+q<n$ we can proceed by induction on $t$. Namely, to prove
the vanishing of $H^p(X,\oqx(t+d))$, $t+d>0$, we may assume that
either $t=0$ or $H^{p+1}(X,\oqqqx(t))=0$

Consider Diagram 1 with exact row and columns.
The horizontal arrows in the diagram are coming
from the cohomology of (2),
so are arrows in the second and fifth (counting from the left) columns
of the diagram (for appropriate $p$, $q$ and $t$).
The other two columns are using arrows from (1).
\par

 From our assumptions concerning cohomology of
$\oqy (t)$ it follows that (remember
that $t+d$ is assumed to be positive)
$$H^p(\oqy (t+d))= H^{p+1}(\oqy (t+d))= 0,$$
and thus
$$H^{p+1}(\oqy (t))\simeq H^p(\oqyx (t+d)).$$
\par

If $t\ne 0$ then also
$H^{p+1}(\oqy (t))=0$ and by inductive
assumption $H^{p+1}(\oqqqx (t))=0$
which yields the desired vanishing of $H^p(\oqx (t+d))$.
\par

If $t=0$ then the map
$$H^p(\oqqqx )\rightarrow H^p(\oqyx (d))\simeq H^{p+1}(\oqy)$$
is just a part of the cupping (3) described in lemma 1.2 so that it is
surjectiv
On the other hand, because of the hard Lefschetz theorem
the cupping map
$$H^{p+1}(\oqqqy) \simeq H^{p+1}(\oqqqx)\rightarrow
H^{p+1}(\oqyx (d))\rightarrow H^{p+2}(\oqy)$$
is injective. This yields the desired vanishing in this case, too.

\bigskip

As a corollary we obtain a special case of a theorem of Flenner
[F, 8.11]:

\proclaim Corollary 1.3. Let $X\subset {\bf P}^N$ be a smooth complete
intersection of dimension $n$.
Then for $0<p<n$, $p+q\neq n$ and $t\in {\bf Z}\setminus \{0\}$ we have
$$H^p(X,\oqx (t))=0.$$

\par

\proclaim Theorem 1.4. Let $X$ and $Y$ be as in the previous theorem $(n \geq
3)$
If $\Omega_Y$ is semistable then $\Omega_X$ is stable.

\par

\noindent{\bf Proof.}
We may assume $det\Omega_Y\equiv{\cal O}(-s)$ and $s\geq d>0$, because
otherwise $X$ is of general type and
the stability of $\Omega_X$ is easily set by Kodaira--Nakano
vanishing. Thus $Y$ is Fano and one may assume that
$PicY\simeq {\bf Z}[\O(1)]$.
Also, we may assume that $s\leq n+2$ and
$(s,d)\ne (n+2,1), (n+1,1)$ since these boundary cases
are clear and related to the hyperplane sections of the projective
space $Y={\bf P}^{n+1}$ and of the smooth quadric $Y={\bf Q}^{n+1}$.
\par
Now we will use a result of Maruyama [Mr, 2.6.1] which asserts that
semistability is preserved under the operation of taking wedge powers.
Therefore, from the semistability of $\Omega_Y$ it follows that
$$H^0(Y,\oqy (t))=0\hbox{ for }t < {q\cdot s\over n+1}.$$
Since $PicX\simeq {\bf Z}[\O(1)]$,
to prove the stability of $\Omega_X$ we have to show
$$H^0(X,\oqx (t))=0\hbox{ for }t\leq {q\cdot (s-d)\over n}.$$
Consider Diagram 2 which is
a  version of the diagram used to prove Theorem 1.1.
We argue as before. Namely, because of our assumptions on
$n$, $s$ and $d$,
$${q\cdot s\over n+1}> {q\cdot (s-d)\over n}$$
and the semistability of $\Omega_Y$ implies
vanishing of $H^0(\oqy (t))$,
for $t$ in question. So we have only to discuss the situation for
$t=0$ or $t=d$ because only then the groups neighbouring
$H^0(\oqx (t))$ are non-zero. But we make exactly the same argument
as before, using cupping, to conclude the appropriate vanishing.
\par

As a corollary we get the following result
for complete intersections in
projective spaces which completes a partial
result obtained by Subramanian in [Sb].

\proclaim Corollary 1.5. Let $X\subset {\bf P}^N$ be a smooth complete
intersection of dimension $n$, $n\geq3$.
Then the cotangent bundle of $X$
is stable.

\bigskip

Now we make similar argument for cyclic coverings of manifolds with
special cohomology.
First we recall some generalities.
\par

Let $\L$ be a line bundle over $Y$
and $D\subset Y$ a smooth divisor
from $\vert \L^{\otimes k}\vert $. By $L$ we denote
the total space $L:=Spec(Sym_Y\L^{\vee})$ of the
bundle ${\L}$ with the projection $\bar\pi: L\lto Y$.
We have a natural coordinate $\xi$ on fibers of $\bar\pi$
coming from the zero section $Y\lto Y_0\subset L$.
Note that $\xi$ is a section of $\bar\pi^*\L$
vanishing along $Y_0$, in other words
$\xi$ descends to $Y$ as the section of
$\bar\pi_*(\bar\pi^*\L)=\L\oplus\O\oplus\L^{\vee}\dots\,$
associated with the $\O$-summand.
Let us assume that the divisor $D$ is defined on a covering
$(U_i)$ by functions $(f_i)$.
Then a cyclic covering associated
to $D$ is a variety $X$ defined
locally in each of $\bar\pi^{-1}(U_i)$ by the
equation $\xi_i^k=f_i$, where $\xi_i$ is the restriction of $\xi$ to
$\bar\pi^{-1}(U_i)$. These equations patch up nicely and define
$X$ globally with a projection $\pi:X\lto Y$.
For more information see [BPV], Section I.17.

(A cyclic covering is a special $k$-section of the line bundle $\L$.
Let us recall that $\pi :X\lto Y$ is a $k$-section of a line
bundle ${\L}$ if $X$ embeds over $Y$ into the total space bundle $L$.)

\proclaim Theorem 1.6.
Let us assume that the pair $(Y,\O(1))$ has special cohomology.
By $\L$ let us denote the line bundle $\O(d)$.
Let $\pi: X\lto Y$
be a $k$-cyclic covering of $X$ branched along
a divisor from $\vert \L^{\otimes k}\vert$.
If $\O_X(1)$ is the pull-back of $\O(1)$ to $X$ then the pair
$(X,\O_X(1)$ has special cohomology.
\par

Proof is very similar to the proof of the Theorem 1.1.
First, we want to understand cupping, now in
$L$. For this purpose we consider
a sequence of normal --- conormal sheaves on $L$:
$$0\lto \bar\pi^*\Omega_Y\lto \Omega_L \lto \Omega^1_{L \vert Y} \otimes
 \bar\pi^*\L^{-1}\lto 0\eqno(4)$$
where $\bar\pi: L\lto Y$ is the projection.
If $(U_i)$ is a covering trivializing both $\L$ and $\Omega Y$
the sequence splits as follows
$$(dy_i)\mapsto (dy_i, d\xi_i)\mapsto (d\xi_i)$$
where $y_i$ are coordinates in $U_i$.

If we twist the above sequence by $\bar\pi^*(\L^t)$, restrict
to $X$ and project to $Y$ the result is
$$\eqalign{
0\lto\Omega_Y\otimes(\L^t\oplus\L^{t-1}\oplus\dots\oplus\L^{t-k+1})
\lto\cr
\pi_*(\Omega_{L\vert X})\otimes\L^t
\lto
\L^{t-1}\oplus\L^{t-2}\oplus\dots\oplus\L^{t-k}
\lto 0}\eqno(5)$$

\proclaim Lemma 1.7. For $t=1\dots k-1$ the split maps and the above sequence
induce a map
$$H^0(Y,\O)\lto
H^1(Y,\Omega_Y\otimes(\L^t\oplus\L^{t-1}\oplus\dots\oplus\L^{t-k+1}))
\lto H^1(Y,\Omega_Y)$$
which defines the first Chern class of $\L$ (that is, the image of the unit $1
\
H^0(Y,{\cal O}_Y )$ is a multiple of the first Chern class of ${\cal L}).$
For $t=k$ the induced map
$$H^0(Y,\O)\lto
H^1(Y,\Omega_Y\otimes
(\L^t\oplus\L^{t-1}\oplus\dots\oplus\L^{t-k+1}))$$
is trivial.

\noindent {\bf Proof of lemma.}
Any section of $\bar\pi^*(\L^t)=\O(tY_0)$ over
$\bar\pi^{-1}(U)$ can be written as
$$\sum_{s\geq -t} g_{is}\xi_i^s$$
where $g_{is}\in \Gamma(U,\L^{-s})$, see [BPV, I.17.2].
Thus, since $(U_i)$ trivializes $\L$,
any section of $\Omega_{L/Y} \otimes \bar\pi^*\L^t$
over $\bar\pi^{-1}(U_i)$ can be written as
$$\sum_{s\geq -t} g_{is}\xi_i^{s}d\xi_i$$
where $g_{is}\in \Gamma(U_i,\O_Y)$.

Over $X\subset L$ we have the following two relations
$$\matrix{\xi_i^k=f_i & \hbox{ and } & k\xi_i^{k-1}d\xi_i=df_i}$$
so that any section of $ \Omega^*_{L/Y}\vert X\otimes\pi^*\L^t$
over $\pi^{-1}(U_i)$ can be written as
$$\sum_{s\geq -t}^{-t+k-2}g_{is}\xi_i^{s}d\xi_i +
g_{is}^0\xi^{-t}df_i.$$
The split map
$$\L^s\lto \pi_*(\Omega_{L/Y}\vert X)
\otimes \L^t \cong \L^{t-1}\oplus\dots
\oplus \L^{t-k+1}$$
for $t-1\leq s\leq t-k+1$
associates to the ``unit'' section of $\L^s$ over $U_i$
the section  $$\xi_i^{-s-1}d\xi_i$$
while for $s=t-k$ the result is
$$\xi_i^{k-t-1}d\xi_i={1\over k}\xi^{-t}df_i.$$

Consequently, the map
$$H^0(Y,\O)\lto
H^1(Y,\Omega_Y\otimes(\L^t\oplus\L^{t-1}\oplus\dots\oplus\L^{t-k+1}))$$
for $t=1\dots k$ can be described as follows
$$\eqalign{1\mapsto ({d\xi_i\over \xi} -{d\xi_j\over \xi_j})=
{1\over k}({df_i\over f_i}-{df_j\over f_j})=
{1\over k}d(log({f_i\over f_j})).}$$
For $t=k$ the above cocycle is trivial as it is the \v Cech boundary of
$$({df_i\over f_i})\in \check C^0((U_i),\Omega_Y\otimes\L^k).$$
For $t\leq k-1$, however, the resulting cocycle is non-trivial and it
defines the first Chern class of $\L$
so that we are done with the lemma.
\medskip

Now we consider an exterior power of the sequence (4),
we restrict it to $X$, project it to $Y$ and twist by a multiple
of $\O(1)$. The result is
$$\eqalign{
0\lto\oqy(t)\oplus\oqy(t-d)\oplus\dots\oplus\oqy(t-(k-1)d)\lto
\pi_*(\Omega^q_{L\vert X})(t)
\cr\lto
\oqqqy(t-d)\oplus\oqqqy(t-2d)\oplus\dots\oplus\oqqqy(t-kd)
\lto 0.}\eqno(6)$$
Arguing as in the proof of Lemma 1.2 we have

\proclaim Lemma 1.8.
For $t=d, 2d,\dots (k-1)d$ the map induced from (6)
$$H^{p-1}(\oqqqy)\lto H^p(\oqy)$$
is the cupping with the first Chern class of $\L$.

Lemma 1.8 implies that some nonzero cohomology groups of $\Omega^q_Y$ do not
con
to the cohomology of $\Omega_{L \vert X}.$ Thus, as a corollary we get the
follo

\proclaim Lemma 1.9. In the situation of Theorem 1.6
for $p+q<dimY$, $p>0$
\item{(i)} $H^p(X,\Omega^q_{L\vert X}(t))=0$ if $t\ne 0, kd.$
\item{(ii)} the induced map
$H^{p}(Y,\pi_*(\Omega^{q}_{L\vert X})(kd))
\lto H^{p}(Y,\Omega_Y^{q-1})$ is injective ( and bijective for $p+q < dim Y
-1$.)
\item{(iii)} the map
$H^p(Y,\oqy)\lto H^p(Y,\pi_*(\Omega^q_{L\vert X}))$
is bijective.
\item{}
\medskip

To conclude Theorem 1.6 we use the following
sequence which is a version of (2) for $X\subset L$
$$0\lto \oqqqx(t)\lto \Omega^q_{L\vert X}(t+kd)
\lto\oqx(t+kd)\lto 0.\eqno(7)$$
We proceed by induction with respect to $t:$ we may assume that the vanishing
is
$\oqqqx(t)$ (for $t$ negative
this is Kodaira---Nakano vanishing) and we want to conclude the vanishing for
 
For $t \ne -kd,0$ and $p+q< dimX, p>0$ we have by (1.9)
$$H^{p+1}(X,\Omega^{q-1}_X (t)) = H^p(X,\Omega^q_{L \vert X}(t+kd)) = 0.$$
For $t = -kd$ we have from 1.9(iii) that $H^p(X,\Omega^q _X) =
H^p(Y;\Omega^q_Y)$
For $t = 0$ the induced map from the cohomology sequence of (7)
$$H^p(X,\Omega_X^{q-1}) \longrightarrow
H^p(X,\Omega^q_{L\vert X}\otimes {\cal O}(X))$$
is a part of cupping (see the proof of (1.1)), so is injective and thus,
because
1.9(ii), it is an isomorphism, therefore $H^p(X,\Omega^q_X(kd)) = 0$.
\bigskip

Now we can prove an extension of Theorem 1.4.
\proclaim Theorem 1.10.
Let $\pi: X\lto Y$ be a cyclic covering which satisfies the assumptions
from Theorem 1.6. If $\Omega_Y$ is semistable and $k\geq 2$
then $\Omega_X$ is stable.

\noindent {\bf Proof}.
Again, as in the proof of (1.4),
we may assume $K_Y=\O(-s)$ for some positive $s$.
Also, we have $dimY+1\geq s\geq (k-1)d> 0$ and
$(s,k,d)\ne (dimY+1,2,1)$
because then we have a double covering of
a projective space by a quadric.
Then
$$K_X=\O_X(-s+(k-1)d)$$
and thus we have to prove
$$H^0(X,\Omega^q_X(t))=0 \hbox{ for }
t\leq {q\cdot(s-(k-1)d)\over dimX}.$$
However, using the semistability of $Y$
(and again Maruyama's result [Mr, 2.6.1]),
and because of (6) we get the vanishing of
$H^0(X,\Omega^q_{L\vert X}(t))$ for
$$t < \hbox{ min }({qs\over dim X}, {(q-1)s\over dimX}+d).$$
On the other hand, in the cohomology of the sequence (7)
(twisted by $\O(-kd)$)
either $H^1(X, \Omega^{q-1}_X(t-kd))$ vanishes or the map
$H^{1}(X,\oqqqx)\lto H^{1}(X,\Omega^q_{L\vert X}\otimes \O(X))$
is injective (as a part of cupping), so that, since
$$t\leq {q\cdot(s-(k-1)d)\over dimX}<
\hbox{ min }({qs\over dim X}, {(q-1)s\over dimX}+d)$$
we have the desired vanishing.
\medskip

\noindent{\bf Remark.}
Note that the ``obvious'' estimate of $H^0(X,\Omega^q_X(t))$
in terms of $H^0(\Omega^{q+1}_{L\vert X}(t+kd))$
which comes from (7) does not give the vanishing in the range
we are interested in.

\beginsection 2. Stability of the tangent bundle
of Fano manifolds with\ $ b_2 = 1$.

Let $X$ always denote
a Fano manifold with $b_2(X) = 1$ and index $r$.
In this section we investigate the problem whether the tangent bundle
$T_X$ is stable, at least for large index.
Let $H$ be the ample generator
of Pic$(X) = {\bf Z}$. For a coherent sheaf $F$  we write for short
$$\mu (F) = {c_1(F)\over rk(F) }$$
instead of $${c_1(F).H^{n-1}}\over rk(F).$$
If the index $ r \geq n $,
then either $ X \simeq {\bf P}^n $ or $ X \simeq
{\bf Q}^n $, the $n$-dimensional quadric, and hence $T_X$ is well-known
to be stable; as a rational homogeneous manifold it even carries a
K\"ahler-Einstein metric.
In order to treat the case of lower index it
is suitable to have a certain condition (ES).
This condition is defined
as follows:

{\narrower\medskip
\item{(ES)}
for $k=1,\dots,r-1$ there exist smooth members
$H_1 ,\dots,H_k \in \vert H \vert $ such
that $H_1 \cap \dots \cap H_k$ is a smooth $(n-k)$-fold
(hence a Fano manifold of index $r-k$).
\medskip}
The condition (ES) is important for
induction purposes and holds for 3-folds
(Sho\-ku\-rov), for 4-folds of index 2 (Wilson, [Wi]) and
generally for index $n-1$ ($=$ {\sl coindex} 2) (Fujita).
\medskip
\proclaim Lemma 2.1.
Assume (ES) for the Fano manifold $X$  and let $r \leq n-1.$
Let  $F \subset T_X $ be a reflexive subsheaf of $T_X$.
Then $\mu(F) < 1$.

\proof Assume $\mu(F) \geq 1$. Let $m = rkF$ , $1\leq m \leq n-1$.
Choose $ H_1 \in \vert H \vert $ smooth and let
$ G_1 = (F_{\vert H_1})^{**}$.
Then we have a morphism
$$\varphi_1 : G_1 \longrightarrow T_{X \vert H_1}  $$
which is generically injective, hence injective
(a priori $ \varphi_1 $
exists only outside codimension 2, then extend).
Via the epimorphism
$$T_{X \vert H_1} \longrightarrow H_{\vert H_1},$$
$\varphi_1 $  induces a map
$\alpha _1 : G_1 \longrightarrow H_{\vert H_1} .$
Now either $\alpha _1 = 0 $,
then we can realize $G_1 $ as a subsheaf of $T_{H_1}$,
or $\alpha_1 \ne 0 $, and then
$c_1(\hbox{Im} \alpha_1 ) = \lambda H\mid H_1 $ for some
$\lambda \leq 1 ,\ $ hence  $\hbox{Ker} \alpha _1 $ is
a subsheaf of $T_{H_1} $
with $c_1 ( \hbox{Ker} \alpha_1 ) \geq c_1(F_{\vert H_1}) - 1 $
(note that $Pic(X) = {\bf Z}$ by Lefschetz,
and consider  Chern classes as numbers).
In this second case we have
$\mu (\hbox{Ker}\alpha _1) \geq 1.$

Let $F_1 = \hbox{Ker} \alpha _1 $
in both cases and repeat the construction
for $F_1 \subset T_{H_1} $.
After $r - 1$ steps we end up with a Fano manifold
$Y = H_1 \cap ...\cap H_{r-1} $ of index 1, $\hbox{ dim} Y = n-r+1$,
and a torsion free sheaf $F' \subset T_Y$.
Observe that at each step we have
$\mu (F_i) \geq 1  $ and $ c_1 (F_i) \leq \hbox{dim} H_i -1$,
hence it follows immediately that always
${rk}F_i\leq dim H_i - 1 $.
Let $k$ be the number of steps with $\alpha _i \neq 0$.
Then $c_1(F') \geq c_1(F) - k $ with
$rk(F') = rk(F) - k $, or $F' = 0$.

First assume $F'\ne 0$. By Reid (2.2 below), $T_Y$ is stable,
hence  we obtain
$${{c_1(F) -k}\over {m-k}}\leq {1\over n-r+1} < 1,$$
hence $c_1(F) < m$, so that $\mu (F) < 1$.

If $F' = 0$, then in particular $ k = m $, thus $r-1 > m $.
Then we stop at the
step where $rk(F_i ) = 1 $.
By a theorem of Wahl, [Wa], (or by [MS, Thm.~8])
we have $c_1 (F_i )\leq 0 $.
Let $k'$  be the number of steps between the first and the $i$-th
where the rank drops, so $k' = m-1 $.
Then $c_1 (F_i ) \geq c_1 (F) -k' $,
hence $c_1 (F)  \leq  m-1 $, proving our claim.
\medskip

For the proof of (2.1) we used

\proclaim Proposition 2.2. {\rm (Reid [Re])} Let X be a Fano manifold
of dimension $n$ with $b_2 (X) = 1 $.
Let $G \subset T_X$ be a proper reflexive subsheaf.
Then $c_1(G) < c_1(X)$ i.e.
$c_1(X)-c_1(G)$ is ample.
In particular, $T_X$ is stable if the index of $X$ is 1.

\proof Indeed, otherwise we would have a subsheaf $G$
of rank $p<n$ with
$c_1(G) \geq c_1(X)$,
hence we obtain a non-zero map $detG \rightarrow \Lambda^pT_X$.
Consequently
$$ 0 \ne H^0(X,\Lambda^p T_X \otimes {det}G^*) =
H^0(X,\Omega^{n-p}_X \otimes {det}G^* \otimes det T_X),$$
which contradicts Kodaira-Nakano vanishing if $c_1(G) > c_1(X)$.
In case of equality $c_1(G) = c_1(X)$,
use also Hodge duality to obtain a contradiction
with vanishing of $H^{n-p}(X,{\cal O}_X)$.
\bigskip

\proclaim Theorem 2.3. Every del Pezzo manifold X
(i.e. X is a Fano n-fold
of index $n-1$) with $b_2=1$ has stable tangent bundle.

\proof Since $\mu (T_X ) = 1 - {1/n}$, $n = dimX $, this follows from
(2.1). Note that (ES) is automatic in this case by Fujita's result.

\medskip \noindent
In particular every Fano 3-fold with $b_2 = 1$
has stable tangent bundle.
Since the only Fano $n$-folds
with index $>n-1$ are projective space and
quadric which have stable tangent bundle we now
turn to Fano n-folds with index $r=n-2$ (i.e.~{\sl coindex} 3).

\proclaim  Proposition 2.4. Let $r=n-2$ and $n \geq 4$.
Assume (ES) and the existence of
a smooth member $H_1 \in \vert H \vert $
such that $T_{H_1}$ is stable. Then $T_X$ is
stable or there exists a reflexive subsheaf
$F \subset T_X $ with $rkF = m$ , $n=2m $
and $\mu (F) = \mu (T_X)$.
Hence $T_X$ is semi-stable and stable for n odd .

\proof
Let $F\subset T_X$ be a proper reflexive subsheaf.
Assume $\mu (F) \geq \mu (T_X) $,
so ${{c_1(F)}/m} \geq {{(n-2)}/n}$.
Since $\mu (F) < 1$ by (2.1), we get $c_1 (F) < m $,
and consequently from both inequalities:
$c_1 (F) = m-1$  and $2m \geq n .$
\smallskip \noindent
Now we have a closer look to the proof of lemma (2.1)
and use the same notations.
If $\alpha _1 : G_1 \rightarrow H_{\vert H_1} $ is 0,
then $G_1 \subset T_{H_1}$,
hence by our assumption we obtain
$$\mu (F) = \mu (G_1) < \mu (T_{H_1}) < \mu (T_X) .$$
So we may assume $ \alpha _1  \not= 0 $.
\smallskip \noindent
Let $F_1 = \hbox{ Ker}\alpha _1$.
Then $\mu (F_1) \geq {{m-2}\over {m-1}}$. On the other hand by
assumption we have
$$\mu (F_1) < {{n-3}\over {n-1}} = \mu (T_{H_1}). $$
Hence $2m < n+1 $, and we conclude $2m = n$,
and $\mu (F) = \mu (T_X)$.

\bigskip \noindent
A direct consequence of 2.4 is

\proclaim Proposition 2.5. Let $X$ be a Fano n-fold of coindex 3
\item{(1)} If $n=4$ and $H^0 (X,\Omega _X^2 (1)) = 0$,
then $T_X$ is stable
\item{(2)} Let  $n=6$ and assume that (ES) holds.
If  $H^0 (X,\Omega _X^3 (2)) = 0$ and there exists
smooth $H_1 \in \vert H \vert$ such that $T_{H_1}$ is stable
then $T_X$ is stable
\item{}

For the proof just observe that (ES) holds for
4-folds of index 2 by Wilson [Wi].

\bigskip \noindent
(2.6) In order to prove {\bf stability} of Fano manifolds of
coindex 3 we make use of
Mukai's classification [Mu].
 Hence from now on we will always assume (ES).
The relevant facts  are collected below. We consider the
Fano 3-fold $Y = H_1 \cap ... \cap H_{n-3} $ which has index 1.
Let $g = {1\over 2} {H^n} + 1 $
be the genus of $X$ (equal to the genus of $Y$).
By the theory of Fano 3-folds, we have $g\leq 12$
and $g\ne 11$. Now $g = 12 $ occurs only in dimension 3,
so we may always assume $g\leq 10$.
The possible values for $g$ are $g = 2 , \dots , 10$.
The linear system $\vert H \vert $ is always
base point free and except for two cases even very ample.
In these exceptional cases, which occur
for g=2 (resp.~g=3), $X$ is a cyclic cover of
$ {\bf P}^n$ (resp.~${\bf Q}^n$), hence by Theorem 1.10,
$T_X$ is stable.
Thus we may assume $\vert H \vert $ very ample. If $ g\leq 5$,
$\vert H\vert $ realises X as a complete intersection
in some projective space, hence by theorem 1.4, $T_X$ is stable.
 So it remains to treat the cases $ 6\leq g \leq 10 $.

\medskip
\proclaim Lemma 2.6. Let $X$ be a Fano manifold of dimension
$n$ and index $n-2$. Let $X_m$, $3\leq m\leq n$,
denote a general linear section of $X$
(that is, an intersection of $n-m$ divisors from the
linear system of $H$, the ample generator of $PicX$).
If $$H^0(X_m ,T_{X_m} ) = 0$$
for $3\leq m\leq {n\over 2}+1$, then $T_X$ is stable.

\proof We may assume that for $m<n$ the bundle $T_{X_m}$ is stable
and $T_X$ is not. Then we
can perform the construction from the proof of 2.1.
In the last part of the proof we considered
$F_i$ which was the "last non-zero sheaf" produced by
 inductive "slicing":
$$ c_1 (F_i)\geq c_1(F)-rk(F)+1.$$
Since by the main theorem of [Wa] (see also [MS])
$c_1 (F_i)\leq 0$, it follows from 2.4
that $c_1(F_i) = 0$ and thus $F_i$
produces a non-zero section in $T_{X_m}$ for some $m\leq n-rk(F)+1$.

\bigskip \noindent
As a direct consequence of 2.6 we obtain
\proclaim Proposition 2.7.  Assume $ n=4 $ and $7\leq g \leq 10 $.
Then $T_X$ is stable.

\proof By [Pr] we have $H^0 (T_Y ) = 0$
for every Fano 3-fold Y with $7 \leq g \leq 10$.
Another proof is easily obtained by applying lemma 2.9.a below.

\proclaim Corollary 2.8. Every Fano $n$-fold of coindex 3 satisfying
(ES) with $n \leq 5$ and
$g \geq 7$ has stable tangent bundle.

\proof Apply (2.4).

\bigskip \noindent
We now deal with all the cases $ 6 \leq g \leq 10$ separately.
First notice that Fano $n$-folds
of coindex 3 and $g=10$  have dimension at most 5,
so this case is already settled.
As to $g=6$ we have:

\proclaim Proposition 2.9. If $g=6$, then $T_X$ is stable.

\proof In this case, the manifold $X$ we are
to deal with (because of (2.4))
can be described as follows. Either
\item{a)} $X$ is a complete intersection
of type (2,1) in $G = G(1,4)$ or
\item{b)} $X$ is a $2:1$ covering of $G$ or of
a 4-dimensional linear section of $G$.
\smallskip \noindent
In case b) the covering is easily seen to be cyclic.
Using the exact sequences (6) and (7) from
Sect.1 and Kodaira-Nakano vanishing we see immediately that
it is sufficient to have
$$ H^0 (Y,\Omega _Y^3 (2)) = 0\ \hbox{ or, respectively, } \
H^0 (Y,\Omega _Y^2 (1)) = 0  $$
for $Y=G=G(1,4)$  or, respectively,  for
$Y$ a 4-dimensional linear section of $G$.
But this follows from the stability of $TY$
since $Y$ is of {\sl coindex} 2.
\smallskip \noindent In case a) we have dim $X= 4$
and by (2.5) it is sufficient to prove
$$ H^0 (X,\Omega _X^2 (1)) = 0 .$$
We have $X$ contained in a hypersurface $Y$ of degree 2 in $G$,
so X is a linear section of Y.
Applying lemma 2.9.a below first to $Y\subset G$ and then
to $X\subset Y$ we obtain
$$H^0(Y,\Omega ^2_Y(1)) = H^0(X,\Omega ^2_X(1)) = 0.$$
The vanishings needed for applying 2.9.a are
$H^1(G,\Omega ^2_G) = 0 = H^1(Y,\Omega ^2_Y)$, which clearly
hold.

\proclaim Lemma 2.9.a.
Let $Y$ be a projective manifold with $Pic(Y) = {\bf Z}$.
Denote by ${\cal O}(1)$ the ample generator.
Let $X$ be a smooth member of the linear system
 $\vert {\cal O}(d) \vert $
for some $d>0$. Fix a number $q <\hbox{dim}Y-1$ such
that $H^1(\Omega^q_Y) = 0$. Then the restriction map
$$ H^0(Y,\Omega^q_Y(c)) \longrightarrow H^0(X,\Omega^q_X(c))$$
is surjective for all $c\leq d$.

\proof We use the following exact sequences (1) and (2) from sect.1.
$$ 0 \longrightarrow \Omega^q_Y(c-d)
\longrightarrow \Omega^q_Y(c) \longrightarrow
\Omega^q_{Y\vert X}(c) \longrightarrow 0$$
and
$$ 0 \longrightarrow \Omega^{q-1}_X(c-d)
 \longrightarrow \Omega^q_{Y\vert X}(c)
 \longrightarrow \Omega^q_X(c) \longrightarrow 0.$$
First let $c<d$. Then it is sufficient to have
$$H^1(X,\Omega^{q-1}_X(c-d)) = 0 \eqno (a_1)$$
$$H^1(Y,\Omega^q_Y(c-d)) = 0 \eqno(a_2)$$
Both claims follow from Kodaira-Nakano vanishing, since $q < $dim$X$.
\par\noindent
If $c=d$, then $(a_2)$ is satisfied by assumption,
however $(a_1)$ has to be replaced by the weaker
statement that the natural map
$$H^1(X,\Omega^{q-1}_X) \longrightarrow
H^1(X,\Omega^q_{Y\vert X}(d)) $$
is injective, which follows from (1.2) and Lefschetz.

\proclaim Corollary 2.10.
Every Fano 4-fold with $b_2 = 1$ has stable tangent bundle.
Every Fano 5-fold with $b_2(X) = 1$ (satisfying (ES)) has stable
tangent bundle except possibly for those of index 2.
\bigskip

We now turn to $g=8$. Then either $X$ is the 8-fold $G(1,5)$
or its linear section.
$G(1,5)$ having stable tangent bundle
as a rational homogeneous manifold (or by (0.1))
it is sufficient by (2.4) to consider
6-dimensional sections .

\proclaim Proposition 2.11. Let $X$ be a 6-dimensional linear
section in  $G = G(1,5)$. Then $T_X$
is stable.

\proof By (2.5) it is sufficient to show $H^0 (\Omega _X^3 (2)) = 0$.
Let $Y$ be a linear section in $G$ and $X$ a linear section of $Y$.
 From the exact sequence
$$ 0 \longrightarrow \Omega_X^3 (2)
\longrightarrow \Omega_{Y\vert X}^4 (3)\longrightarrow \Omega_X^4 (3)
\longrightarrow 0 \eqno(S)$$
we deduce that it is sufficient to prove
$$ H^0(X,\Omega_{Y\vert X}^4(3)) = 0 \eqno(1)$$
This in turn follows from $$ H^0(Y,\Omega_Y^4(3)) = 0 \eqno(2)$$
$$H^1(Y,\Omega_Y^4(2)) = 0. \eqno(3)$$
We first prove (2). Since $\Omega_Y^4(3)$ is
a subbundle of $\Omega_{G\vert Y}^5(4)$,
it is sufficient to prove $$H^0(Y,\Omega_{G\vert Y}^5(4)) = 0.$$
This in turn is guaranteed by
$$H^0(G,\Omega_G^5(4)) = H^1(G,\Omega_G^5(3)) = 0.$$
In order to prove (3) we first observe that tensoring (S)
by ${\cal O}(-1)$
and substituting $X$ by $Y$,
$Y$ by $G$, and taking cohomology, it suffices to show
$$H^1(Y,\Omega_{G\vert Y}^4(2)) = 0 \eqno(4)$$
$$H^2(Y,\Omega_Y^3(1)) = 0 \eqno(5)$$
Now (4) is immediate since from (0.1) it follows
$$H^1(G,\Omega^4_G(2)) = H^2(G,\Omega_G^4(1)) = 0.$$
To prove (5) we use an analogous sequence to (S) which yields
an exact sequence
$$H^2(Y,\Omega_Y^2)  \buildrel \beta \over \longrightarrow
H^2(Y,\Omega_{G\vert Y}^3(1)) \buildrel \alpha \over
\longrightarrow
H^2(Y,\Omega_Y^3(1)) \longrightarrow H^3(Y,\Omega_Y^2).$$
Observe $H^3(\Omega^2_Y) = 0$ since $b_5(Y) = 0$.
Moreover $h^{2,2}(Y) = h^{2,2}(G) = 2$
and since $$H^2(G,\Omega_G^3(1)) = H^3(G,\Omega_G^3(1)) = 0,$$
we have
 $$H^2(G,\Omega_{G\vert Y}^3(1)) = H^3(G,\Omega_G^3) = {\bf C}^2.$$
The composite map
$$H^2(G,\Omega_G^2) \longrightarrow H^2(Y,\Omega^2_Y)
\buildrel \beta \over \longrightarrow H^2(Y,\Omega_{G\vert Y}^3(1))
\longrightarrow H^3(G,\Omega_G^3)$$
is cupping by the fundamental class of $Y \subset G$
(Lemma 1.2), hence an isomorphism ($h^{2,2}= h^{3,3}$).
This proves surjectivity of $\beta$, hence $\alpha = 0$, hence (5).
\bigskip

If now $g=9$ , $X$ is a rational homogeneous 6-fold,
thus $T_X$ is stable. The last case is then
$g=7$. Here $X$ is a 10-dimensional homogeneous manifold,
a so-called spinor variety $S_{10} $
(see [Mu]), hence $T_X$ is stable,
or $X$ is a linear section of $S_{10}$ . The stability of $T_X$ of a
6-dimensional section (hence also of a 7-dimensional section by (2.4))
follows from

\proclaim Lemma 2.12. Let $X_k$ be a smooth $k$-dimensional linear section of
th
10-dimensional spinor variety. If $k\geq 5$ then
$$H^0(X_k,\Omega^3_{X_k}(2))=H^1(X_k,\Omega^3_{X_k}(1))=
H^1(X_k,\Omega^2_{X_k}(1))=0.$$

\par

Proof of the lemma is by descending induction: for $k=10$ this follows
from Snow's computations. Therefore we may assume all vanishing to
hold for $k+1$. Also, we will use the knowledge of the following
Betti numbers of $X_{10}$ and thus, by Lefschetz, of $X_{k+1}$:
$b_2=b_4=1$, $b_3=b_5=0$.
\par

To prove the last vanishing of the lemma we use an exact sequence of sheaves
(see the sequence (2) in the proof of theorem 1.1)
$$0\lto\Omega_{X_k}\lto\Omega^2_{X_{k+1}\vert X_k}(1)
\lto\Omega^2_{X_k}(1)\lto 0$$
and the resulting exact sequence of cohomology
$$H^1(X_k,\Omega_{X_k})\lto H^1(X_k,\Omega^2_{X_{k+1}\vert X_k}(1))
\lto H^1(X_k,\Omega^2_{X_k}(1))\lto
H^2(X_k,\Omega_{X_k}).$$
The last term in the above sequence is 0 since $b_3(X_k)=0$.
The first arrow in the sequence is surjective because it is a part of cupping
$$H^1(X_{k+1},\Omega_{X_{k+1}})\lto
H^1(X_k,\Omega_{X_k})\lto H^1(X_k,\Omega^2_{X_{k+1}\vert X_k}(1))\lto
H^2(X_{k+1},\Omega^2_{X_{k+1}})$$ (see lemma 1.2),
the cupping is isomorphic since $b_2(X_{k+1})=b_4(X_{k+1})=1$, moreover
the restriction
$H^1(X_{k+1},\Omega_{X_{k+1}})\lto H^1(X_k,\Omega_{X_k})$
is an isomorphism by Lefschetz and the kernel of
$H^1(X_k,\Omega^2_{X_{k+1}\vert X_k}(1))\lto
H^2(X_{k+1},\Omega^2_{X_{k+1}})$ is $H^1(X_{k+1},\Omega^2_{X_{k+1}}(1))$
(see diagram 1) so it vanishes by our inductive assumption.
This proves vanishing of $H^1(X_k,\Omega^2_{X_k}(1))$.
\par

To prove vanishing of $H^1(X_k,\Omega^3_{X_k}(1))$ a similar
exact sequence of cohomology is applied:
$$H^1(X_k,\Omega^3_{X_{k+1}\vert X_k}(1))
\lto H^1(X_k,\Omega^3_{X_k}(1))\lto
H^2(X_k,\Omega^2_{X_k})\lto H^2(X_k,\Omega^3_{X_{k+1}\vert X_k}(1)).$$
The last arrow in the above sequence is injective as a part of cupping
(see 1.2) while the first term above can be bounded by the cohomology
of the following exact sequence (see the sequence (1) from the proof of 1.1):
$$0\lto\Omega^3_{X_{k+1}}\lto\Omega^3_{X_{k+1}}(1)\lto
\Omega^3_{X_{k+1}\vert X_k}(1)\lto 0.$$
Indeed, by the inductive assumption $H^1(X_{k+1},\Omega^3_{X_{k+1}}(1))=0$
and because
$$b_5(X_{k+1})=H^2(X_{k+1},\Omega^3_{X_{k+1}})=0$$
it follows that $H^1(X_{k+1},\Omega^3_{X_{k+1}\vert X_k}(1))=0$.
Therefore $H^1(X_k,\Omega^3_{X_k}(1))=0$.
\par

Now the vanishing of $H^0(X_k,\Omega^3_{X_k}(2))$ is easy: we have the
following two exact sequences of cohomology (cf.~Diagram 2):
$$H^0(X_{k+1},\Omega^3_{X_{k+1}}(2))\lto
H^0(X_{k+1},\Omega^3_{X_{k+1}\vert X_k}(2))\lto
H^1(X_{k+1},\Omega^3_{X_{k+1}}(1))$$
$$H^0(X_{k},\Omega^3_{X_{k+1}\vert X_k}(2))\lto
H^0(X_{k},\Omega^3_{X_{k}}(2))\lto
H^1(X_{k},\Omega^2_{X_{k}}(1))$$
and the vanishing follows.

\medskip It remains to prove the stability of the tangent bundle of an
8-dimensional linear section $X$. From (2.4) we see that for this purpose
it is sufficient to prove
\medskip
\proclaim Lemma 2.13
 \it Let $X$ be an 8-dimensional linear section of the 10-dimensional spinor
manifold $S$. Then $H^0(\Omega_X ^4(3)) = 0.$
 \medskip \noindent \bf Proof.\rm Choose a smooth linear section $Y\subset S$
su
For both inclusions we will use several times the exact sequences (1)
and (2) from sect.1 , with $d = 1.$  We will quote these sequences just
by (1) and (2).
By (2) with $q=4,t=3$ it is sufficient to show $$ H^0(\Omega^5_Y(4) \vert X) =
0
$$
Using (1) this comes down to show $$ H^0(\Omega^5_Y(4)) = 0 \eqno(b) $$
$$H^1(\Omega^5_Y(3)) = 0 \eqno(c)$$
Using again (2) , now for $Y \subset X$, (b) comes down to
$$ H^0(\Omega^6_S(5) \vert Y) = 0.$$
This is justified by the vanishings (using (1))
$$H^0(\Omega^6_S(5)) = H^1((\Omega^6_S(4)) = 0$$
which both hold by [Sn2,3].
In order to prove (c) we use (2) with $q=5,t=3$ and we have to show
$$H^1(\Omega^6_S(4)\vert Y) = 0 \eqno(d)$$
$$H^0(\Omega^6_Y(4)) = 0 \eqno(e) $$
Now (d) follows by (1) from the vanishings
$$H^1(\Omega^6_S(4)) = H^2(\Omega^6_S(3)) = 0$$
and (e) is verified by $$H^0(\Omega^7_S(5) \vert Y)) = 0$$
which in turn is guaranteed by
$$H^0(\Omega^7_S(5)) = H^1(\Omega_S^7(4)) = 0.$$
All the needed vanishings on $S$ follow again from [Sn2,3].

\medskip
In summary we obtain
\medskip
\proclaim Theorem 2.14.
Let X be a n-dimensional Fano manifold of coindex 3, with $b_2 = 1$.
Assume (ES). Then $T_X$ is stable.

\beginsection References.

\item{[BPV]} Barth, W., Peters, C., van de Ven, A.:
Compact complex surfaces. Springer 1984
\item{[F]} Flenner, H.: Divisorenklassengruppen
quasihomogener Singularit\"aten,
Crelle's J. 328, 128---160 (1981)
\item{[H]} Hartshorne, R.: Algebraic Geometry. Springer 1977
\item{[La]} Lascoux, A.: Cohomologie de la grassmannienne.
Preprint 1987
\item{[Mr]} Maruyama, M.: The theorem of Grauert-M\"ulich-Spindler.
Math.~Ann.~255, 317-333 (1981)
\item{[Mt]} Matsumura, H.: Geometric structure on the cohomology
rings in abstract algebraic geometry.
Mem.~Coll.~Sci.~Univ.~ Kyoto 32, 33-84 (1959)
\item{[MS]} Mori, S., Sumihiro, N.: On Hartshorne's conjecture.
J.Math. Kyoto Univ.18, 523-533 (1978)
\item{[Mu]} Mukai, S.: New classification of Fano threefolds
and Fano manifolds of coindex 3.
Proc.~Natl.~Acad.~Sci.~ USA 86, 3000-3002, (1989).
\item{[Pr]} Prokhorov, Y.: Groups of automorphisms of Fano manifolds.
Uspiechy Mat.~Nauk 45 (1990)
\item{[Re]} Reid, M.: Bogomolov's theorem $c_1^2 \leq 4c_2$. in:
Intl.~Symp.~on Algebraic Geometry, Kyoto 1977, 623-642 , Kinokuniya
\item{[S]} Sommese, A.: Hyperplane sections of projective surfaces.
Duke Math~J. 46, 377-401 (1979)
\item{[Sn]}  Snow, D.: Cohomology of twisted
holomorphic forms on Grassmann manifolds and
quadric hypersurfaces. Math.~Ann 276, 159-176 (1986)
\item{[Sn2]} Snow,D.: Vanishing theorems on compact hermitian symmetric
spaces. Math.Z.198,1-20 (1988)
\item{[Sn3]} Snow,D.: Letter to the authors , aug. 12, 1993
\item{[St]} Steffens, A.: Stabilit\"at des Tangentialb\"undels
dreidimensionaler Fanomannigfaltig\-kei\-ten. Thesis, 1993
\item{[Sb]} Subramanian, Stability and existence of K\"ahler---Einstein
metric, Math.~Ann.~291 (1991), 573---577.
\item{[Ti]} Tian, G.: On the stability of the tangent
bundles of Fano varieties. Intl.~J.~Math. 3, 401-413 (1992)
\item{[Wa]} Wahl, J.: A cohomological characterization of ${\bf P}^n$
Inv.~Math.~72, 315-32 (1983)
\item{[We]} Weyman, J.: The equations of conjugacy
classes of nilpotent matrices. Inv.~Math. 98, 229-245 (1989)
\item{[Wi]} Wilson, P.M.H.: Fano fourfolds of
index greater than one. Crelle's J.379, 172-181 (1987)

\vskip 1.5 cm
{\centerline{{\sl Mathematisches Institut,
Universit\"at Bayreuth, 95440 Bayreuth, Germany}}}
{\centerline{{\sl Instytut Matematyki,
Uniwersytet Warszawski, 02-097 Warszawa, Poland}}}

\end



\def\O{{\cal O}}
\def\L{{\cal L}}
\def\oqy{\Omega^q_Y}

\def\oqqqy{\Omega^{q-1}_Y}
\def\oqx{\Omega^q_X}
\def\oqqx{\Omega^{q+1}_X}
\def\oqqqx{\Omega^{q-1}_X}
\def\oqyx{\Omega^q_{Y\vert X}}
\def\oqqyx{\Omega^{q+1}_{Y\vert X}}

\def\oqqqyx{\Omega^{q-1}_{Y\vert X}}
\def\Diag{\def\normalbaselines{\baselineskip25pt
\lineskip3pt\lineskiplimit3pt}
        \matrix}
\def\lto{\rightarrow}

\hsize=7.6in \vsize=8.5in
\hoffset=-.6in
\nopagenumbers

\midinsert
$$
\Diag{
H^p(\Omega^{q-2}_X(t-d))&&           &    &               &
&H^{p+1}(\Omega^{
\downarrow             &&           &    &               &    &\downarrow
  H^p(\oqqqyx(t))&&H^p(\oqy(t+d))   & =  & 0\hfill       &
&H^{p+1}(\oqqqyx(t
  \downarrow     && \downarrow      &    &               &    & \downarrow
H^p(\oqqqx(t))&\lto&H^p(\oqyx (t+d))&\lto&H^p(\oqx (t+d))&\lto&H^{p+1}(\oqqqx
(t
  \downarrow     && \downarrow      &    &               &    & \downarrow
H^{p+1}(\Omega^{q-2}_X(t+d))&&H^{p+1}(\oqy (t))&  &       &
&H^{p+2}(\Omega^{
                 && \downarrow      &    &               &    &
                 &&H^{p+1}(\oqy (t+d))& =& 0\hfill       &    &
}
$$
\centerline{Diagram 1}\vglue.1in
\endinsert

\midinsert
$$
\Diag{
H^0(\Omega^{q-2}_X(t-2d))&&          &    &               &
&H^1(\Omega^{q-2}
\downarrow       &&                 &    &               &    &\downarrow
  H^0(\oqqqyx(t-d))&&H^0(\oqy(t))   &    &               &
&H^1(\oqqqyx(t-d))
  \downarrow     && \downarrow      &    &               &    & \downarrow
H^0(\oqqqx(t-d))&\lto&H^0(\oqyx (t))&\lto&H^0(\oqx (t))&\lto&H^1(\oqqqx (t-d))
  \downarrow     && \downarrow      &    &               &    &
H^1(\Omega^{q-2}_X(t-2d))&&H^1(\oqy (t-d))&  &            &    &
}
$$
\centerline{Diagram 2}\vglue.1in
\endinsert

\end